\documentclass[leqno]{amsart}
       
\usepackage{amsthm,amsmath,amsfonts,amssymb}

\def\R{{\mathbb R}}
\def\N{{\mathbb N}}

\def\Sch{{\mathcal S}} 
\def\virgp{\raise 2pt\hbox{,}}

\def\bu{{\bf u}}

\def\({\left(}
\def\){\right)}
\def\<{\left\langle}
\def\>{\right\rangle}
\def\le{\leqslant}
\def\ge{\geqslant}

\def\Tend#1#2{\mathop{\longrightarrow}\limits_{#1\rightarrow#2}}

\def\d{{\partial}}
\def\a{{\tt a}}

\def\eps{\varepsilon}
\def\h{h}

\def\om{\omega}
\def\si{{\sigma}}

\def\O{\mathcal O}

\DeclareMathOperator{\RE}{Re}
\DeclareMathOperator{\IM}{Im}
\DeclareMathOperator{\Op}{\mathrm{Op}}
\DeclareMathOperator{\DIV}{\mathrm{div}}

\def\defn{\mathrel{:=}}

\theoremstyle{plain}
\newtheorem{theorem}{Theorem}[section]
\newtheorem{lemma}[theorem]{Lemma}
\newtheorem{corollary}[theorem]{Corollary}

\theoremstyle{definition}
\newtheorem{definition}[theorem]{Definition}

\theoremstyle{remark}
\newtheorem{remark}[theorem]{Remark}

\numberwithin{equation}{section}

\title[Loss of regularity for NLS]{Loss of regularity for
 supercritical nonlinear Schr\"odinger equations}
 \author[T. Alazard]{Thomas Alazard}
\address{CNRS \& Universit\'e
  Paris-Sud\\ Math\'ematiques\\ B\^at. 425\\ 91405 
  Orsay cedex\\ France}
\email{Thomas.Alazard@math.cnrs.fr}
\author[R. Carles]{R{\'e}mi Carles}
\address{CNRS \& Universit\'e Montpellier~2\\ Math\'ematiques\\ CC 051\\ 
   Place Eug\`ene Bataillon\\ 34095
  Montpellier cedex 5\\ France}
\email{Remi.Carles@math.cnrs.fr}

\begin{document}

\begin{abstract}
We consider the nonlinear Schr\"odinger equation with defocusing,
smooth, 
nonlinearity. Below the critical Sobolev regularity, it is known that
the Cauchy problem is ill-posed. We show that this is even worse,
namely that there is a loss of regularity, in the spirit of the
result due to 
G.~Lebeau in the case of the wave equation. As a consequence, the
Cauchy problem for energy-supercritical equations is not well-posed
in the sense of Hadamard. 
We  reduce the problem to a
supercritical WKB analysis. For super-cubic, smooth nonlinearity,
this analysis is new, and relies on the introduction of a modulated
energy functional \emph{\`a la} Brenier. 
\end{abstract}

\thanks{Support by the ANR project SCASEN is acknowledged.} 
\maketitle

\section{Introduction}
\label{sec:intro}

We consider the following defocusing nonlinear Schr\"odinger equation
on $\R^n$: 
\begin{equation}
  \label{eq:nls}
  i\d_t \psi+\frac{1}{2}\Delta \psi = |\psi|^{2\si}\psi\quad ;\quad \psi_{\mid
  t=0}=\varphi, 
\end{equation}
where $\si\ge 1$ is an integer, so that the nonlinearity is
smooth. It is well-known that the critical Sobolev regularity
corresponds to the value given by scaling arguments,
\begin{equation*}
  s_c\defn \frac{n}{2} - \frac{1}{\si}\cdot 
\end{equation*}
Throughout this paper, we assume $s_c>0$. If $\varphi\in H^s(\R^n)$ with
$s \ge s_c$, then the Cauchy problem \eqref{eq:nls} is locally
well-posed in $H^s(\R^n)$ \cite{CW90}. On the other hand, if
$s<s_c$, then the Cauchy problem \eqref{eq:nls} is ill-posed
\cite{CCT2} (see also the appendices in \cite{BGTENS,CaARMA}). 
The
worst phenomenon proved in \cite{CCT2} is the norm inflation. For
$0<s<s_c$, one can 
find a sequence $(\psi^\h)_{0<\h\le 1}$ of solutions to
\eqref{eq:nls} and $0<t^\h\to 0$, such that $\varphi^\h \in
\Sch(\R^n)$  and  
\begin{equation*}
  \|\varphi^\h\|_{H^s}\Tend \h 0 0 \quad ;\quad
  \|\psi^\h(t^\h)\|_{H^s}\Tend \h 0 +\infty.
\end{equation*}
In this paper, we prove the stronger result:
\begin{theorem}\label{theo:new}
Let $\si\ge 1$. Assume that $s_c=n/2-1/\si >0$, and let
$0<s<s_c$. 
There exists a family $(\varphi^\h)_{0<\h \le 
  1}$ in ${\mathcal S}({\mathbb R}^n)$ with 
\begin{equation*}
  \|\varphi^\h\|_{H^{s}({\mathbb R}^n)} \to 0 \text{ as
  }\h \to 0, 
\end{equation*}
a solution $\psi^\h$ to
\eqref{eq:nls} and $0<t^\h \to 0$, such that: 
\begin{equation*}
  \|\psi^\h(t^\h)\|_{H^{k}({\mathbb R}^n)} \to +\infty \text{ as }\h \to
 0\, , \ \forall k> \frac{s}{1+\si(s_c-s)}\cdot
\end{equation*}
\end{theorem}
Note that this result is not bound to the case $x\in \R^n$: for
instance, it remains valid on a compact manifold, see
\S\ref{sec:geometrie}.
\smallbreak

In the case $\si =1$ and $n\ge 3$, this result was established in
  \cite{CaARMA}. It followed from a supercritical WKB
  analysis for the cubic nonlinear Schr\"odinger equation, which had
  been justified by E.~Grenier~\cite{Grenier98}. For $\si \ge 2$,
  adapting the results of \cite{Grenier98} seems to be a much more
  delicate issue, and a rigorous analysis in this setting for $n\le 3$
  has been given very recently \cite{AC-BKW}. An important remark, in
  the proof of Theorem~\ref{theo:new} that we present here, is that it
  is not necessary to justify WKB analysis as precisely as in
  \cite{Grenier98}, or \cite{AC-BKW}, to obtain this
  loss of regularity. From this point of view, our proof is very
  simple. On the other hand, it can be considered as highly
  nonlinear: it relies on a \emph{quasilinear} analysis, as opposed
  to the \emph{semilinear} analysis in \cite{CCT2} (see also
  Remark~\ref{rem:semiquasi} at the end of this paper). In the opinion of
  the authors, the proof of Theorem~\ref{theo:new} is at least as
  interesting as the result itself. 
\begin{remark}
 Shortly after this work was completed, an alternative proof was given
 by L.~Thomann \cite{ThomannAnalytic}, based on the justification of
 WKB analysis in an 
 analytic setting. This approach allows to consider focusing
 nonlinearities (the nonlinearity is treated as a semilinear
 perturbation in spaces based on analytic regularity), unlike the
 method followed in the present paper. On 
 the other hand, the virial identity shows that for supercritical
 focusing nonlinearities, blow-up can happen for arbitrary small data
 in $H^s$ ($s<s_c$) and arbitrary small times
 (see e.g. \cite[Exercise~3.63]{TaoDisp}).       
\end{remark}
\smallbreak

This result is to be compared with the main result in
\cite{Lebeau05}, which we recall with notations adapted to make the
comparison with the Schr\"odinger case easier. For supercritical wave
equations  
\begin{equation*}
  \(\d_t^2-\Delta\)u + u^{2\si+1}=0,
\end{equation*}
G.~Lebeau shows
that one can find a \emph{fixed} initial datum in $H^s$, and a
sequence of times $0<t^\h\to 0$, such that the $H^k$ norms of the
solution are unbounded along the sequence $t^\h$, for $k\in
]I(s),s]$. The expression for  $I(s)$ is related to the critical
Sobolev exponent
\begin{equation*}
  s_{\text{sob}}=\frac{n}{2}\frac{\si}{\si+1}\virgp
\end{equation*}
which corresponds to the embedding
$H^{s_{\text{sob}}}(\R^n)\subset L^{2\si +2}(\R^n)$. In
\cite{Lebeau05}, we find:
\begin{equation}\label{eq:Ileb}
  I(s)=1\text{ if } 1<s\le s_{\text{sob}}\quad ;\quad 
I(s)=\frac{s}{1+\si(s_c-s)}\text{ if } s_{\text{sob}}\le s<s_c.
\end{equation}
Note that we have
\begin{equation}
  \label{eq:sature}
\frac{s_{\text{sob}}}{1+\si(s_c-s_{\text{sob}})}=1.
\end{equation}
The approach in \cite{Lebeau05} consists in using an \emph{anisotropic}
scaling, as opposed to the isotropic scaling used in
\cite{Lebeau01,CCT2}. Compare Theorem~\ref{theo:new} with
the approach of \cite{Lebeau05}. Recall that
\eqref{eq:nls} has two important (formally) conserved quantities: mass
and energy,
\begin{equation}
  \label{eq:conserv}
  \begin{aligned}
&M(t) = \int_{\R^n}|\psi(t,x)|^2dx\equiv M(0),\\
&E(\psi(t)) = \frac{1}{2}\int_{\R^n}|\nabla \psi(t,x)|^2dx
+\frac{1}{\si+1}\int_{\R^n}|\psi(t,x)|^{2\si+2}dx\equiv E(\varphi).
\end{aligned}
\end{equation}
In view of \eqref{eq:sature}, we obtain, for
$H^1$-supercritical nonlinearities: 
\begin{corollary}\label{cor:energy}
Let $n\ge 3$ and $\si >\frac{2}{n-2}$. 
There exists  $(\varphi^\h)_{0<\h \le 
  1}$ in ${\mathcal S}({\mathbb R}^n)$ with 
\begin{equation*}
\| \varphi^\h \|_{H^1}+ \| \varphi^\h\|_{L^{2\si+2}}\rightarrow 0
\text{ as
  }\h \to 0, 
\end{equation*}
a solution $\psi^\h$ to
\eqref{eq:nls} and $0<t^\h \to 0$, such that: 
\begin{equation*}
  \|\psi^\h(t^\h)\|_{H^{k}({\mathbb R}^n)} \to +\infty \text{ as }\h \to
 0\, , \ \forall k>1.
\end{equation*}
\end{corollary}
We thus get the analogue
of the result of G.~Lebeau when $I(s)=1$, with the drawback that we
consider a \emph{sequence} of initial data only. The information that
we don't have for Schr\"odinger equations, and which is available for
wave equations, is the finite speed of propagation, that is used in
\cite{Lebeau05} to construct a fixed initial datum; see also the
discussion in \S\ref{sec:geometrie}.  On the other hand,
our approach involves an \emph{isotropic} scaling, which is recalled
and generalized in Section~\ref{sec:reduc}. Moreover, our range for $k$ is
broader when $1<s<s_{\rm sob}$, and also, we allow the range $0<s\le
1$, for which no analogous result is available for the wave
equation. Note that  
unlike in \cite{Lebeau05}, we perform no linearization in our
analysis (the properties of the analogous linearized operator are not
as  interesting
in the case of Schr\"odinger equations): despite the fact that
for fixed $\eps$, \eqref{eq:nls} is a 
semilinear equation, we consider a quasilinear system to
prove our main result. 
\smallbreak

Before going further into details, let us focus on the notion of solution
to \eqref{eq:nls}. In view of Theorem~\ref{theo:new}, we assume that
the initial data are in the 
Schwartz class: $\varphi\in \Sch(\R^n)$. Then \eqref{eq:nls} has a
local smooth solution: for all $s>n/2$, there exists $T_s>0$ such
that \eqref{eq:nls} has a unique solution $\psi \in
C([-T_s,T_s];H^s)$.  If $n\le 2$, then 
\eqref{eq:nls} has a global smooth solution, $\psi \in C(\R;H^s)$ for
all $s\ge 0$, and the identities \eqref{eq:conserv} hold for all
time. The same is 
true when $n=3$ and $\si =1$. These results are established in
\cite{GV84}. In the $H^1$-critical three dimensional case ($n=3$ and
$\si=2$), it is proved in \cite{CKSTTAnnals} that solutions with $H^s$
regularity ($s>1$) remain in $H^s$ for all time; the same is true in
the four dimensional case ($n=4$ and $\si=1$), from \cite{RV}. 
On the other hand, if the nonlinearity is
$H^1$-supercritical ($\si > \frac{2}{n-2}$), then it is
not known in general whether the solution remains smooth for all time
or not. In Theorem~\ref{theo:new}, for $\si=1$, the solution $\psi^h$
is a smooth 
solution, that remains smooth up to time $t^h$, thanks to a result due
to E.~Grenier \cite{Grenier98}.  \emph{A priori}, the solution we
consider in Theorem~\ref{theo:new} is a weak solution:
\begin{definition}\label{def:weak}
  Let $\varphi \in H^1\cap L^{2\si+2}(\R^n)$. A (global) weak solution to
  \eqref{eq:nls}  is a function $\psi \in C(\R;L^2)\cap
  L^\infty(\R; H^1\cap L^{2\si+2})\cap C_w (\R; H^1\cap L^{2\si+2})$
  solving \eqref{eq:nls} in 
  ${\mathcal D}'(\R\times \R^n)$, and such that:
  \begin{itemize}
  \item $\|\psi(t)\|_{L^2}= \|\varphi\|_{L^2}$, $\forall t\in \R$.
\item $E(\psi(t))\le E(\varphi)$, $\forall t\in \R$.
  \end{itemize}
\end{definition}
From \cite{GV85c}, for $\varphi \in \Sch(\R^n)$, \eqref{eq:nls} has a
global weak solution. The proof in \cite{GV85c} is based on Galerkin
method. We use a different construction, as in \cite{Lebeau05}, which
is described in Section~\ref{sec:og}. Note that when the
nonlinearity is $H^1$-subcritical, then the weak solution is unique,
and coincides with the strong solution. Recall also that the existence
of blowing-up solutions in the $H^1$-supercritical case is open so
far. On the other hand, if the nonlinearity is focusing, many
results are available (see \cite{TaoDisp} for an overview of the
subject, and similar problems for other dispersive equations).  
\smallbreak

Note that in view of
Definition~\ref{def:weak}, Corollary~\ref{cor:energy} is sharp. 
\smallbreak

As in \cite{CaARMA}, the idea for the proof of Theorem~\ref{theo:new}
consists in reducing the analysis to a supercritical WKB analysis,
for an equation of the form:
\begin{equation}
  \label{eq:nlssemi}
  i\eps \d_t u^\eps +\frac{\eps^2}{2}\Delta u^\eps =
  |u^\eps|^{2\si}u^\eps\quad ;\quad u^\eps(0,x) = a_0(x). 
\end{equation}
The parameter $\eps$ goes to zero. The above equation is
supercritical as far as geometrical optics is concerned: if one plugs
an approximate solution of the form
\begin{equation*}
  v^\eps \sim e^{i\phi/\eps}\(\a_0+\eps \a_1+\eps^2 \a_2+\ldots\)
\end{equation*}
into the equation, then closing the systems of equations for
$\phi,\a_0,\a_1,\ldots$ is a very delicate issue (see
e.g. \cite{CaBKW}). In the case 
$\si=1$, this issue was resolved by E.~Grenier
\cite{Grenier98}. However, the argument in \cite{Grenier98}
relies very strongly on the fact that the nonlinearity is defocusing,
and cubic at the origin. In \cite{AC-BKW}, we have
proposed an approach that justifies  WKB analysis in Sobolev spaces for
\eqref{eq:nlssemi} for any $\si\ge 2$, in space dimension $n\le 3$
(higher dimensions could also be considered with the same proof, up to
considering sufficiently large values of $\si$). In
\cite{PGX93} and \cite{ThomannAnalytic}, WKB analysis was justified in spaces
based on analytic regularity, in a periodical setting and for analytic
manifolds, respectively. As noticed in
\cite{AC-BKW}, analytic regularity is necessary to justify WKB
analysis with a focusing nonlinearity. 
The analytic regularity essentially allows to
view the nonlinearity as a semilinear perturbation, and to construct
an approximate solution that solves \eqref{eq:nlssemi} up to a source
term of order $e^{-\delta/\eps}$ for some $\delta>0$. Yet,
such a justification is not needed to prove
Theorem~\ref{theo:new}; see \S\ref{sec:reduc}. In this paper, we use
a functional 
that yields sufficiently many informations to infer
Theorem~\ref{theo:new}. This functional may be viewed as a
generalization of the one used in \cite{LinZhang} in the cubic case,
following an idea introduced by Y.~Brenier \cite{BrenierCPDE}. The
general form for this modulated energy functional was announced in
\cite{LinZhang}. However, we will see that making the corresponding
analysis rigorous is not straightforward, since we consider weak
solutions. 
\smallbreak

The main idea of the proof of Theorem~\ref{theo:new} consists in
noticing that for an $\eps$-independent initial datum $a_0$ in
\eqref{eq:nlssemi}, the solution $u^\eps$ becomes $\eps$-oscillatory
for times of order $\O(1)$ as $\eps \to 0$. This phenomenon is typical
of supercritical r\'egimes, as far as geometrical optics is concerned
(see also \cite{CG05}). This crucial step is stated in
Theorem~\ref{theo:reduc}, which in turn is proved thanks to the above
mentioned modulated energy functional.  

\smallbreak
We end this introduction with a remark concerning the study of the 
Cauchy problem for \eqref{eq:nls}. 
From \cite{CCT2}, 
it is known that the Cauchy problem is not well posed in $H^s(\R^n)$
for $0<s<s_{c}$.  
Yet, one can try to solve the Cauchy problem by searching the
solutions in a larger space. Denote $H^\infty = \cap_{s>0}
H^s(\R^n)$. Recall the notion of well-posedness in the sense of Hadamard:
\begin{definition}\label{def:WP}
Let $s\ge k\ge 0$. The Cauchy problem for \eqref{eq:nls} 
is well posed from $H^{s}(\R^n)$ to $H^{k}(\R^n)$ if, for all bounded
subset $B\subset H^{s}(\R^n)$, there exist  
$T>0$ and a Banach space $X_T\hookrightarrow C([0,T];H^{k}(\R^n))$
such that: \\
(1) For all $\varphi\in B\cap H^\infty $, \eqref{eq:nls} has a unique
  solution $\psi \in
  C([0,T];H^\infty)$. \\
(2) The mapping $\varphi\in (H^\infty,\| \cdot \|_{B})
  \mapsto \psi\in X_T$ is continuous. 
\end{definition}
The following result is a direct consequence of our analysis (see
Remark~\ref{rem:hadamard}).
\begin{corollary}\label{cor:CauchyEnd}
Let $n\ge 1$ and $\si\ge 1$ be such that $\sigma>2/n$. 
The Cauchy problem for \eqref{eq:nls} is not well posed from
$H^s(\R^n)$ to $H^k(\R^n)$ for all $(s,k)$ such that 
$$
0<s<s_c=\frac{n}{2}-\frac{1}{\sigma}\virgp \qquad k>
\frac{s}{1+\si(s_c-s)}\cdot $$
\end{corollary}

\section{Reduction of the problem}
\label{sec:reduc}

Let $0<s<s_c$ and $a_0\in \Sch(\R^n)$. For a sequence $\h$ aimed at
going to zero, 
consider the family of initial data
\begin{equation}\label{eqs:varphi}
  \varphi^\h(x) =\h^{s-\frac{n}{2}}a_0\(\frac{x}{\h}\). 
\end{equation}
Let $\eps = \h^{\si(s_c-s)}$. By assumption, $\eps$ and $\h$ go
simultaneously to zero. Define the function $u^\eps$ by the relation:
\begin{equation}\label{eqs:ueps}
  u^\eps (t,x) = \h^{\frac{n}{2}-s} \psi^\h \( \h^2 \eps t,\h x\). 
\end{equation}
Then \eqref{eq:nls} is equivalent to \eqref{eq:nlssemi}.
Note that we have the relation:
\begin{equation*}
  \|\psi^\h(t)\|_{\dot H^m} = \h^{s-m}\left\|u^\eps\(
  \frac{t}{\h^2\eps}\)\right\|_{ \dot H^m}.
\end{equation*}
Our aim is to show that for some $\tau>0$ independent of $\eps$, 
\begin{equation}\label{eq:epsosc}
 \liminf_{\eps \to 0}\eps^k \left\|u^\eps\(
  \tau\)\right\|_{ \dot H^k}>0,\quad \forall k\ge 0. 
\end{equation}
Back to $\psi$, this will yield $t^\h= \tau \h^2\eps$ and 
\begin{equation*}
  \|\psi^\h(t^\h)\|_{\dot H^k} \gtrsim \h^{s-k}\eps^{-k}=
  \h^{s-k(1+\si(s_c-s))}.
\end{equation*}
To complete the above reduction, note that in view of
Theorem~\ref{theo:new}, we only have to prove \eqref{eq:epsosc} for
$k\in ]0,1]$. Indeed, for $k>1$, there exists $C_k>0$ such that 
\begin{equation*}
  \|f\|_{\dot H^1}\le C_k \|f\|_{L^2}^{1-1/k}\|f\|_{\dot
  H^k}^{1/k},\quad \forall f\in H^k(\R^n).
\end{equation*}
This inequality is straightforward thanks to Fourier analysis. Note
also that thanks to the conservation of mass for $u^\eps$, we have:
\begin{equation*}
  \|u^\eps(t)\|_{\dot H^1}\le C_k \|a_0\|_{L^2}^{1-1/k}\|u^\eps(t)\|_{\dot
  H^k}^{1/k}.
\end{equation*}

Up to replacing $a_0$ with $|\log \h|^{-1}a_0$ (for instance), the
analysis of this section shows that Theorem~\ref{theo:new} follows
from:
\begin{theorem}\label{theo:reduc}
Let $n\ge 1$, $a_0\in \Sch(\R^n)$ be non-trivial, and $\si\ge 1$.
There exists a solution
$u^\eps$ to \eqref{eq:nlssemi} and $\tau>0$ such that for all $k\in ]0,1]$,
\begin{equation*}
  \liminf_{\eps \to 0} \left\lVert \lvert\eps D_{x}\rvert^k u^\eps\(
  \tau\)\right\rVert_{L^2}>0,\quad \text{where }D_x =-i\nabla.
\end{equation*}
\end{theorem}
\begin{remark}\label{rem:hadamard}
As we will see, the previous conclusion holds for all family of smooth 
solutions $u^\eps$ defined on a time interval independent of
$\eps$. In particular,  
Corollary~\ref{cor:CauchyEnd} also 
follows from this analysis. To see this, suppose, by
contradiction, that the Cauchy problem is well  
posed from $H^{s}(\R^n)$ to $H^{k}(\R^n)$. 
Since the family of initial data given by \eqref{eqs:varphi} 
is bounded in $H^s(\R^n)$, the first point in Definition~\ref{def:WP}
implies that  
the solutions $\psi^\h$ are defined for a time interval $[0,T]$
independent of~$h$. As a result, the function  
$u^\eps$, as given by \eqref{eqs:ueps}, is defined for $t\in
[0,T/(\eps h^2)]$ with value in $H^\infty(\R^n)$, and hence on the
fixed  time interval $[0,T]$. Then, Theorem~\ref{theo:og} implies that
there exists $\tau>0$ such that  
$\liminf \| |\eps D_{x}|^{k}u^\eps(\tau)|\|_{L^2}>0$. Back to
$\psi^\h$ this yields the existence of a sequence $\tau^h$ such that  
$\| \psi^h (\tau^h)\|_{H^k}$ tends to $+\infty$, which contradicts the
continuity given by the second point of  
the definition.
\end{remark}

\begin{remark}\label{rema:Iunif}
  If we could prove Theorem~\ref{theo:reduc} for $k=1$ only, then back
  to $\psi^h$, this would yield Theorem~\ref{theo:new} for $I(s)<k\le s$,
  where $I(s)$ is given by \eqref{eq:Ileb}, like in
  \cite{Lebeau05}. 
\end{remark}
Consider the case $k=1$, and recall that the conservation of energy for
$u^\eps$ reads, as long as $u^\eps$ is a strong solution of
\eqref{eq:nlssemi}:  
\begin{equation*}
  E^\eps(t) = \frac{1}{2}\int_{\R^n}|\eps\nabla u^\eps(t,x)|^2dx
+\frac{1}{\si+1}\int_{\R^n}|u^\eps(t,x)|^{2\si+2}dx\equiv E^\eps(0).
\end{equation*}
At time $t=0$, the first term (kinetic energy) is of order
$\O(\eps^2)$, while the second (potential energy) is dominating, of
order $\O(1)$. Therefore, the game consists in showing that there
exists $\tau>0$, time at which the kinetic energy is of the order of
the total (initial) energy as $\eps\to 0$. 
\smallbreak

Some important features of the proof of this result 
can be revealed by analyzing the linear case with variable
coefficients: 
\begin{equation*}
  i\eps \d_t u^\eps +\frac{\eps^2}{2}\Delta u^\eps =
  V(x) u^\eps\quad ;\quad u^\eps(0,x) = a_0(x). 
\end{equation*}
Introduce the operator $H^\eps\defn -(\eps^2 / 2)\Delta + V(x)$, 
so that $u^\eps (t)=e^{-i t H^\eps/\eps}a_0$. 
Now, let $\Op_\eps (q)$ be a semiclassical pseudo-differential operator 
with symbol $q(x,\xi)\in S^1_{1,0}$. 
Since $e^{i t H^\eps/\eps}$ is unitary, by means of Egorov's Theorem
(see \cite{Martinez}), we obtain 
\begin{align*}
\| \Op_\eps (q)u^\eps\|_{L^2}
=\| e^{i t H^\eps/\eps}\Op_\eps (q)e^{-i t H^\eps/\eps} a_0 \|_{L^2}=
\|\Op_\eps (q\circ \Phi_t)a_{0}\|_{L^2}+\mathcal{O}(\eps),
\end{align*}
where $\Phi_{t}$ is the Hamiltonian flow associated with $H^\eps$. 
For small times, one can
relate $\Phi_t$ to the solution $\phi(t,x)$ of the Hamilton--Jacobi equation:
\begin{equation*}
\partial_t \phi +\frac{1}{2}|\nabla\phi|^2+V(x)=0\quad;\quad \phi(0,x)=0,
\end{equation*}
by the identity $\Phi_t(x,\xi)=(X(t,x)+t\xi\, , \,
\xi+(\nabla\phi)(t,X(t,x))),$ 
where $X$ satisfies 
\begin{equation*}
  \partial_{t}X(t,x)=(\nabla\phi)(t,X(t,x))\quad ;\quad X(0,x)=x.
\end{equation*}
Hence, with $q(x,\xi)=i\xi$, we infer
$$
\| \eps\nabla u^\eps(t)\|_{L^2} = \| (\nabla\phi)(t,X(t,x))
a_0\|_{L^2}+\O(\eps), 
$$
so that the kinetic energy is of order $\O(1)$ provided that
$(\nabla\phi)(t,X(t,\cdot)) a_0\neq 0$.   

\smallbreak

The previous argument can be made explicit for  
the harmonic oscillator 
\begin{equation}
  \label{eq:nlslinear}
  i\eps \d_t u^\eps_\ell +\frac{\eps^2}{2}\Delta u^\eps_\ell =
 \frac{|x|^2}{2} u^\eps_\ell\quad ;\quad u^\eps_\ell(0,x) = a_0(x). 
\end{equation}
\begin{lemma}\label{lem:linear}
Let $n\ge 1$, and $a_0\in \Sch(\R^n)$ (non-trivial). There exists
$\tau>0$ such that 
the solution $u^\eps_\ell$ to 
\eqref{eq:nlslinear} satisfies
\begin{equation*}
  \liminf_{\eps \to 0} \left\lVert \eps \nabla u^\eps_\ell\(
  \tau\)\right\rVert_{L^2}>0.
\end{equation*}
\end{lemma}
\begin{proof}
The standard WKB approach yields, at leading order, the following
approximate solution:
\begin{equation*}
  v^\eps_\ell(t,x) = a_\ell(t,x)e^{i\phi_\ell(t,x)/\eps},
\end{equation*}
where $\phi_\ell$ and $a_\ell$ are given by an eikonal equation and a
transport equation. 
Since we consider an harmonic
oscillator, we can compute $\phi_\ell$ and $a_\ell$ explicitly:
\begin{align*}
  &\d_t \phi_\ell +\frac{1}{2}|\nabla \phi_\ell|^2
  +\frac{|x|^2}{2}=0 ;\  \phi_{\ell \mid t=0}=0:
\quad \phi_\ell(t,x)
  = \frac{-|x|^2}{2}\tan t.\\ 
  &\d_t a_\ell +\nabla \phi_\ell\cdot \nabla a_\ell +\frac{1}{2}a_\ell
  \Delta \phi_\ell=0  ;\  a_{\ell \mid t=0}=a_0:
 \ a_\ell(t,x)=
  \frac{1}{(\cos 
  t)^{n/2}}a_0\(\frac{x}{\cos t}\). 
\end{align*}
Energy estimates then yield (see for instance \cite[\S 3]{CaBKW} for more
details): 
\begin{equation*}
  \|\eps\nabla u^\eps_\ell -\eps\nabla
  v^\eps_\ell\|_{L^\infty([0,T];L^2)} \le C_T \eps,\quad \forall T\in
  \left[ 0,\frac{\pi}{2}\right[. 
\end{equation*}
Since 
\begin{equation*}
 \liminf_{\eps \to 0} \left\lVert \eps \nabla v^\eps_\ell\(
  t\)\right\rVert_{L^2}= \sin t \left\lVert x a_0
  \right\rVert_{L^2}, \quad \forall t\in
  \left[ 0,\frac{\pi}{2}\right[, 
\end{equation*}
the lemma follows easily. 
\end{proof}
The strategy for proving Theorem~\ref{theo:reduc} is the same: we
compare with the limit system. For nonlinear  
Schr\"odinger equation, the eikonal equation which
gives the phase is coupled  
to the transport equation: the limiting system reads
\begin{equation}
  \label{eq:limiteS}
  \left\{
    \begin{aligned}
     \d_t \phi +\frac{1}{2}|\nabla \phi|^2+ |a|^{2\si}=0
\quad &;\quad \phi_{\mid t =0}=0,\\
\d_t a +\nabla\phi\cdot \nabla a +\frac{1}{2}a
\Delta \phi = 0\quad   &;\quad a_{\mid
  t=0}=a_0. 
 \end{aligned}
\right.
\end{equation}
By introducing $v=\nabla\phi$, one can transform this system into a
quasilinear system of nonlinear equations.  
An important feature of the system thus obtained is that it does not 
enter into the classical framework of symmetric hyperbolic systems
for $\si\ge 2$.  
Nevertheless, one can solve the Cauchy problem~\eqref{eq:limiteS} for
all $\si\ge 1$  
by a nonlinear change of variable. This is done in
\S\ref{sec:muk} following  
an idea due to T.~Makino, S.~Ukai and S.~Kawashima~\cite{MUK86}. 
\smallbreak

For the general case $\si\ge 1$, we establish a modulated energy
estimate, following the 
pioneering  
work of Y.~Brenier \cite{BrenierCPDE}. The idea consists in obtaining
an estimate  
for the $L^{2}$ norm of $\eps\nabla a^\eps$ where $a^\eps$ is the
modulated unknown function $a^\eps \defn u^\eps e^{-i\phi/\eps}$.  
It is found that 
$a^\eps$ satisfies
\begin{align*}
i\eps \Bigl( \partial_t a^\eps +\nabla\phi\cdot\nabla a^\eps
+\frac{1}{2}a^\eps \Delta\phi\Bigr) 
+ \frac{\eps^2}{2}\Delta a^\eps
= \( \left| a^\eps \right|^{2\si} - \left|a\right|^{2\si}\) a^\eps .
\end{align*}
For $\si=1$, one can obtain estimates uniform in $\eps$, that is 
$$
\left\| \eps\nabla a^\eps\right\|_{L^\infty([0,T];L^2)} 
+ \left\| |a^\eps|^2-|a|^2 \right\|_{L^\infty([0,T];L^2)}
=\O(\eps),
$$
by an integration by parts argument. Again, this is based on the
hyperbolicity in the case $\si=1$ (see \cite{LinZhang} for an
application  
of this idea to the {G}ross--{P}itaevskii equations). 
Using a modulated energy functional adapted to our problem, we prove
the estimate (see
Theorem~\ref{theo:og} below): 
$$
\left\| \eps\nabla a^\eps\right\|_{L^\infty([0,T];L^2)} 
+ \left\| \(|a^\eps|^2-|a|^2\)\(|a^\eps|^{\si-1}
  + |a|^{\si -1}\)  \right\|_{L^\infty([0,T];L^2)}
=\O(\eps).
$$
This is enough to prove Theorem~\ref{theo:reduc} for $k=1$. Note that
this suffices to infer
Corollary~\ref{cor:energy}.  Finally,
to cover the range $k\in ]0,1]$, we microlocalize the previous
estimate  by means of wave packets operator.

\section{The limiting system}
\label{sec:muk}

Being optimistic, one would try to mimic the approach of E.~Grenier
\cite{Grenier98}, and write the solution $u^\eps$ to
\eqref{eq:nlssemi} as 
$  u^\eps = a^\eps e^{i\phi^\eps /\eps}$, where
\begin{equation}
  \label{eq:ideal}
  \left\{
    \begin{aligned}
     \d_t \phi^\eps +\frac{1}{2}|\nabla \phi^\eps|^2+ |a^\eps|^{2\si}=0
\quad &;\quad \phi^\eps_{\mid t =0}=0,\\
\d_t a^\eps +\nabla\phi^\eps\cdot \nabla a^\eps +\frac{1}{2}a^\eps
\Delta \phi^\eps = i\frac{\eps}{2}\Delta a^\eps\quad   &;\quad a^\eps_{\mid
  t=0}=a_0. 
 \end{aligned}
\right.
\end{equation}
Considering the unknown $v^\eps = \nabla \phi^\eps$ instead of
$\phi^\eps$, the first step in the analysis would be to solve
\begin{equation}
  \label{eq:idealv}
  \left\{
    \begin{aligned}
     \d_t v^\eps +v^\eps\cdot \nabla v^\eps+ \nabla\(|a^\eps|^{2\si}\)=0
\quad &;\quad v^\eps_{\mid t =0}=0,\\
\d_t a^\eps +v^\eps\cdot \nabla a^\eps +\frac{1}{2}a^\eps
\DIV v^\eps = i\frac{\eps}{2}\Delta a^\eps\quad   &;\quad a^\eps_{\mid
  t=0}=a_0. 
 \end{aligned}
\right.
\end{equation}
In \cite{Grenier98}, E.~Grenier considers the unknown $\bu^\eps =
(v^\eps,\RE a^\eps,\IM a^\eps)\in \R^{n+2}$. It solves a partial
differential equation of the form
\begin{equation*}
  \d_t \bu^\eps +\sum_{j=1}^n
  A_j(\bu^\eps)\d_j \bu^\eps
  = \frac{\eps}{2} L  
  \bu^\eps .
\end{equation*}
In the case $\si=1$, the left-hand side of
the above equation defines a symmetric quasilinear hyperbolic system in
the sense of Friedrichs, with a constant symmetrizer. The linear
operator $L$ corresponds to the term $i\Delta$ on the right hand side
of \eqref{eq:idealv}: it is skew-symmetric, and does not appear in the
energy estimates. Therefore, one can construct a smooth solution to
\eqref{eq:idealv} on some time interval $[0,T]$ with $T>0$ independent
of $\eps$. In the case $\si\ge 2$,  the symmetrizer of
\cite{Grenier98} would become
\begin{equation*}
  S=\left(
    \begin{array}[l]{cc}
      \frac{1}{4\si |a^\eps|^{2\si -2}}I_n&0\\
 0& I_2
    \end{array}
\right).
\end{equation*}
For $a^\eps \in L^2(\R^n)$, this matrix is not uniformly bounded, and
this is why the analysis in \cite{Grenier98} is restricted to
nonlinearities which are defocusing, and cubic at the origin.
\smallbreak

This apparent lack of hyperbolicity is not a real problem for the
homogeneous nonlinearity that we consider, provided that we analyze
the limiting system only:
\begin{equation}
  \label{eq:limite}
  \left\{
    \begin{aligned}
     \d_t \phi +\frac{1}{2}|\nabla \phi|^2+ |a|^{2\si}=0
\quad &;\quad \phi_{\mid t =0}=0,\\
\d_t a +\nabla\phi\cdot \nabla a +\frac{1}{2}a
\Delta \phi = 0\quad   &;\quad a_{\mid
  t=0}=a_0. 
 \end{aligned}
\right.
\end{equation}
The above restriction remains apparently valid for this
system: in the presence of vacuum (zeroes of $a$), the symmetrizer $S$
is singular. This may lead to a loss of regularity in the energy
estimates. However, we shall see that thanks to the special structure of
\eqref{eq:limite}, we can construct solutions to \eqref{eq:limite} in
Sobolev spaces of sufficiently large order.  
Following an idea due to T.~Makino, S.~Ukai and S.~Kawashima
\cite{MUK86}, we prove: 
\begin{lemma}\label{lem:muk}
  Let $a_0\in \Sch(\R^n)$. There exists $T>0$ such that
  \eqref{eq:limite} has a unique solution $(\phi,a)\in
  C^\infty([0,T]\times \R^n)^2$, with $(\phi,a)\in
  C([0,T],H^s)^2$ for all $s\ge 0$. Moreover, $ \<x\>^s\nabla \phi \in
  C([0,T],L^2)$ for all $s\ge 0$, where $\<x\>=(1+|x|^2)^{1/2}$.   
\end{lemma}
\begin{proof} 
 Differentiating the first equation in \eqref{eq:limite}, we first
 consider:
\begin{equation}
  \label{eq:limitev}
  \left\{
    \begin{aligned}
     \d_t v +v\cdot \nabla v+ \nabla\(|a|^{2\si}\)=0
\quad &;\quad v_{\mid t =0}=0,\\
\d_t a +v\cdot \nabla a +\frac{1}{2}a
\DIV v = 0\quad   &;\quad a_{\mid
  t=0}=a_0. 
 \end{aligned}
\right.
\end{equation}
Adapting the idea of \cite{MUK86}, consider the unknown
$(v,u)=(v,a^\si)$. Even though the map $a\mapsto a^\si$ is not
bijective, this will suffice
to prove the lemma. The pair $(v,u)$  solves:
\begin{equation}
  \label{eq:limitev+}
  \left\{
    \begin{aligned}
     \d_t v +v\cdot \nabla v+ \nabla\(|u|^{2}\)=0
\quad &;\quad v_{\mid t =0}=0,\\
\d_t u +v\cdot \nabla u +\frac{\si}{2}u
\DIV v = 0\quad   &;\quad u_{\mid
  t=0}=a_0^\si\in \Sch(\R^n). 
 \end{aligned}
\right.
\end{equation}
This system is hyperbolic symmetric, with a constant symmetrizer. To
see this, denote $\bu=(v,\RE u,\IM u)^T\in
\R^{n+2}$. Equation~\eqref{eq:limitev+} is of the form  
\begin{equation*}
\d_t \bu +\sum_{j=1}^n A_j(\bu)\d_j \bu =0,
\end{equation*}
where the matrices $A_j \in {\mathcal M}_{n+2}(\R)$ are such that
$SA_j$ are symmetric, for 
\begin{equation*}
S=\(
\begin{array}{cc}
 \si {\rm I}_n   &  0  \\
   0  &  4     {\rm I}_2
\end{array}
\)\in {\mathcal S}_{n+2}(\R).
\end{equation*}
From classical theory on hyperbolic symmetric quasilinear systems (see
  e.g. \cite{AlinhacGerardUS,Taylor3}), 
there exist  $T>0$ and a unique solution $(v,u)\in C^\infty([0,T]\times
\R^n)^2$, which is in $C([0,T],H^s)^2$ for all $s\ge 0$.  
  The fact that $\<x\>^s v\in
  C([0,T],L^2)$ follows easily by considering the momenta of $u$ and
  $v$. Now that $v$ is known,
  we define 
  $a$ as the solution of the transport equation
  \begin{equation*}
   \d_t a +v\cdot \nabla a +\frac{1}{2}a
\DIV v = 0\quad   ;\quad a_{\mid
  t=0}=a_0.  
  \end{equation*}
The function $a$ has the regularity announced in
Lemma~\ref{lem:muk}. We check that $a^\si$ solves the second equation
in \eqref{eq:limitev+}. Since $v$ is a smooth coefficient, by
uniqueness for this linear equation, we have
$u=a^\si$. Therefore, $(v,a)$ solves \eqref{eq:limitev}. To conclude, we
notice that $v$ is irrotational, so there exists $\widetilde \phi$
such that $v=\nabla \widetilde \phi$. Setting $\phi=\widetilde \phi +
F$, where $F=F(t)$ is a function of time only, $(\phi,a)$ solves
\eqref{eq:limite}. Uniqueness follows from the uniqueness for
\eqref{eq:limitev+}. 
\end{proof}
\begin{remark}
  The proof shows that if we assume only $a_0\in H^s(\R^n)$ with
  $s>n/2+1$, then $u,v \in C([0,T];H^s)$. We infer $a \in
  C([0,T];H^{s-1})$: the possible loss of regularity due to the lack
  of hyperbolicity for \eqref{eq:limite} remains limited.  
\end{remark}
\begin{remark}
The nonlinear change of unknown function, $u=a^\si$, suggests that the
above approach cannot be adapted to study \eqref{eq:idealv}, since we
have to deal with the term $i\Delta a^\eps$, and prevent the loss of
regularity that it may cause in the energy estimates. 
\end{remark}

\section{Semi-classical limit}
\label{sec:og}
Introduce the hydrodynamic variables:
  \begin{equation*}
    \rho = |a|^2\quad ;\quad \rho^\eps = |u^\eps|^2 \quad ;\quad
    J^\eps = \IM\(\eps \overline u^\eps \nabla u^\eps\). 
  \end{equation*}
The main result of this section is:
\begin{theorem}\label{theo:og}
  Let $n\ge1$, and $\si\ge 1$ be an integer. Let $(v,a)\in
  C([0,T];H^\infty)^2$ given by Lemma~\ref{lem:muk}, where $v=\nabla
  \phi$. Then we have the following estimate:
  \begin{align*}
    \left\| (\eps\nabla -iv)u^\eps\right\|_{L^\infty([0,T];L^2)}^2 
+ \left\| \(\rho^\eps-\rho\)^2\((\rho^\eps)^{\si-1}+\rho^{\si-1} \)
    \right\|_{L^\infty([0,T];L^1)} 
=\O(\eps^2).  
  \end{align*}
\end{theorem}
Note that the above quantities are well-defined for weak
solutions. We outline the argument in a formal proof, which is then made
rigorous. 

\begin{proof}[Formal proof]
For $y\ge 0$, denote 
\begin{align*}
  &f(y)=y^\si\quad ;\quad F(y)=\int_0^y f(z)dz = \frac{1}{\si +1}y^{\si
  +1}\quad ;\\
& G(y)= \int_0^y zf'(z)dz= yf(y)-F(y) =
  \frac{\si}{\si+1}y^{\si +1}.  
\end{align*}
We check that
$(\rho^\eps,J^\eps)$ satisfies, for $\si\ge 1$:
\begin{equation}
  \label{eq:hydroq}
  \left\{
    \begin{aligned}
      &\d_t \rho^\eps + \DIV J^\eps =0.\\
&\d_t J^\eps_j+\frac{\eps^2}{4}\sum_k \d_k\(4\RE \d_j \overline u^\eps
\d_k u^\eps - \d_{jk}^2 \rho^\eps\) +\d_j G 
(\rho^\eps) =0. 
    \end{aligned}
\right.
\end{equation}
As suggested in \cite[Remark~1, (2)]{LinZhang}, introduce the
modulated energy functional:
\begin{equation*}
  H^\eps(t) = \frac{1}{2}\int_{\R^n} \left| (\eps \nabla
  -iv)u^\eps\right|^2dx + \int_{\R^n} \(F(\rho^\eps) -
  F(\rho)-(\rho^\eps -\rho)f(\rho)\)dx. 
\end{equation*}
Denote 
\begin{equation*}
  K^\eps(t)= \frac{1}{2}\int_{\R^n} \left| (\eps \nabla
  -iv)u^\eps\right|^2dx . 
\end{equation*}
Integrations by parts, which are studied in more detail below,  yield:
\begin{equation*}
  \frac{d }{dt}H^\eps(t) = \O\(K^\eps +\eps^2\) -\int_{\R^n}
  \(G(\rho^\eps)-G(\rho) -(\rho^\eps-\rho)G'(\rho)\) \DIV v\,  dx. 
\end{equation*}
We check that there exists $K>0$ such that 
\begin{equation*}
  \left\lvert G(\rho^\eps)-G(\rho) -(\rho^\eps-\rho)G'(\rho)\right
  \rvert\le K \left\lvert F(\rho^\eps)-F(\rho) -(\rho^\eps-\rho)F'(\rho)\right
  \rvert . 
\end{equation*}
Therefore,
\begin{equation*}
  H^\eps (t ) \le H^\eps (0 ) + C \int_0^t \(
  H^\eps (s )+\eps^2\)ds . 
\end{equation*}
We infer by Gronwall lemma that $H^\eps (t ) =\O(\eps^2)$
so long as it is defined.
Finally, Taylor's formula yields, since $F''(y)=f'(y)=\si y^{\si-1}$:
\begin{equation*}
  F(\rho^\eps)-F(\rho) -(\rho^\eps-\rho)F'(\rho) =
  \si \(\rho^\eps-\rho\)^2\int_0^1 (1-\theta)\( \rho
  +\theta\(\rho^\eps-\rho\) \)^{\si-1}d\theta.
\end{equation*}
The estimate, for all $\theta\in [0,1]$,
\begin{equation*}
 \( \rho
  +\theta\(\rho^\eps-\rho\) \)^{\si-1} = \( (1-\theta)\rho
  +\theta\rho^\eps \)^{\si-1} \ge (1-\theta)^{\si-1}\rho^{\si-1} +
  \theta^{\si-1} \(\rho^\eps\)^{\si-1}
\end{equation*}
shows that
there exists $c>0$ such that 
\begin{equation}\label{eq:convexfinal}
 H^\eps (t )\ge K^\eps(t) + c\int_{\R^n} (\rho^\eps-\rho)^2 \(
  (\rho^\eps)^{\si-1} + \rho^{\si-1}\)dx.
 \end{equation}
The result of Theorem~\ref{theo:og} follows. 
\end{proof}

\begin{proof}[Rigorous proof] 
  In general, the above integrations by parts do not make sense for all
  $t\in [0,T]$, since we consider weak solutions only. Note however
  that for $\si \ge 2$ and $n\le 3$, the analysis in \cite{AC-BKW}
  shows that we can work with strong solutions, so the following
  analysis is not needed in this case (for $\si=1$ and $n\ge 1$, the
  same holds true, from \cite{Grenier98}). Also, if one is just
  interested in proving Corollary~\ref{cor:CauchyEnd} by
  contradiction, no further justification is needed for integrations
  by parts, and one can skip the end of this section. 

We work on a sequence of global strong solutions,
  converging to a weak solution. For $(\delta_m)_m$ a sequence of positive
  numbers going to zero,  introduce the saturated
  nonlinearity, defined for $y\ge 0$: 
  \begin{equation*}
    f_m(y)= \frac{y^\si}{1+\(\delta_m y\)^{\si}}\cdot 
  \end{equation*}
Note that $f_m$ is a symbol of degree $0$. 
For fixed $m$, we have a global strong solution $u^\eps_m\in
C(\R;H^1)$ to:
\begin{equation}
  \label{eq:nlssemireg}
  i\eps \d_t u^\eps_m +\frac{\eps^2}{2}\Delta u^\eps_m =
  f_m\(|u^\eps_m|^2\) u^\eps_m\quad ;\quad u^\eps_m(0,x) = a_0(x). 
\end{equation}
As $m\to \infty$, the sequence $(u^\eps_m)_m$ converges to a weak
solution of \eqref{eq:nlssemi} (see \cite{GV85c,Lebeau05}). 
For $y\ge 0$, introduce also
\begin{align*}
  F_m(y)=\int_0^y f_m(z)dz\quad ;\quad
 G_m(y)= \int_0^y zf'_m(z)dz= yf_m(y)-F_m(y) .  
\end{align*}
The mass and energy associated to $u_m^\eps$ are conserved: 
\begin{align*}
  &M_m^\eps(t)=\int |u^\eps_m(t,x)|^2dx \equiv \|a_0\|_{L^2}^2.\\
&E^\eps_m(t) = \frac{1}{2}\|\eps\nabla u^\eps_m(t)\|_{L^2}^2 +
\int_{\R^n} F_m \( |u^\eps_m(t,x) |^2\) dx \equiv E^\eps_m(0). 
\end{align*}
Moreover, the solution is in $H^2(\R^n)$ for all time: $u^\eps_m\in
C(\R;H^2)$. To see this, we use an idea due to T.~Kato
\cite{Kato87,KatoNLS}, and consider $\d_t u^\eps_m$. Energy estimates
show that $\d_t u^\eps_m \in C(\R;L^2)$, since $f_m$ is a symbol of
degree $0$. Using
\eqref{eq:nlssemireg} and the boundedness of $f_m$, we
infer $\Delta u^\eps_m \in C(\R;L^2)$. 
\smallbreak

We consider the hydrodynamic variables:
  \begin{equation*}
    \rho^\eps_m = |u^\eps_m|^2 \quad ;\quad
    J^\eps_m = \IM\(\eps \overline u^\eps_m \nabla u^\eps_m\). 
  \end{equation*}
From the above discussion, we have:
\begin{equation}
  \label{eq:reghydrom}
  \rho^\eps_m (t)\in
W^{2,1}(\R^n)\text{ and }J^\eps_m(t)\in 
W^{1,1}(\R^n),\quad \forall t\in \R.
\end{equation}
The analogue of \eqref{eq:hydroq} is:
\begin{equation}
  \label{eq:hydroqm}
  \left\{
    \begin{aligned}
      &\d_t \rho^\eps_m + \DIV J^\eps_m =0.\\
&\d_t (J^\eps_m)_j+\frac{\eps^2}{4}\sum_k \d_k\(4\RE \d_j \overline u^\eps_m
\d_k u^\eps_m - \d_{jk}^2 \rho^\eps_m\) +\d_j G_m
(\rho^\eps_m) =0. 
    \end{aligned}
\right.
\end{equation}
Introduce the modulated energy functional ``adapted to
\eqref{eq:nlssemireg}'':
\begin{equation*}
  H^\eps_m(t) = \frac{1}{2}\int_{\R^n} \left| (\eps \nabla
  -iv)u^\eps_m\right|^2dx + \int_{\R^n} \(F_m(\rho^\eps_m) -
  F_m(\rho)-(\rho^\eps_m -\rho)f_m(\rho)\)dx. 
\end{equation*}
Notice that this functional is not exactly adapted to
\eqref{eq:nlssemireg}, since the limiting quantities (as $\eps \to 0$)
$\rho$ and $v$ are constructed with the nonlinearity $f$ and not the
nonlinearity $f_m$. We also distinguish the kinetic part:
\begin{equation*}
  K^\eps_m(t) = \frac{1}{2}\int_{\R^n} \left| (\eps \nabla
  -iv)u^\eps_m\right|^2dx. 
\end{equation*}
Thanks to the conservation of energy for $u^\eps_m$, we have:
\begin{align*}
  \frac{d}{dt}K^\eps_m =& -\frac{d}{dt} \int F_m(\rho^\eps_m)dx
  +\frac{1}{2}\int |v|^2 \d_t \rho^\eps_m +\int \rho^\eps_m v \cdot
  \d_t v\\
&-\int J^\eps_m \cdot \d_t v - \int v\cdot \d_t J^\eps_m. 
\end{align*}
Using Lemma~\ref{lem:muk}, \eqref{eq:reghydrom} and
\eqref{eq:hydroqm}, (licit) integrations by parts yield:
\begin{align*}
  \frac{d}{dt}K^\eps_m =& -\frac{d}{dt} \int F_m(\rho^\eps_m)dx
  -\frac{1}{2}\int |v|^2 \DIV J^\eps_m - \sum_{j,k}\int \rho^\eps_m v_j v_k
  \d_j v_k \\
&- \int \rho^\eps_m \nabla f(\rho)\cdot v 
+ \int (v\cdot \nabla v)\cdot J^\eps_m + \int  \nabla f(\rho)\cdot
  J^\eps_m\\
 -\sum_{j,k} \int& \d_k v_j \RE \( \eps \d_j \overline
  u^\eps_m \eps \d_k u^\eps_m \)
 -\frac{\eps^2}{4}\int \nabla\(\DIV
  v\)\cdot \nabla \rho^\eps_m +\int \rho^\eps_m v\cdot \nabla
  f_m(\rho^\eps_m). 
\end{align*}
Proceeding as in \cite{LinZhang}, we have:
\begin{align*}
  \eps^2\int \nabla \(\DIV v\)\cdot \nabla \rho^\eps_m &=  \eps\int
  \nabla \(\DIV v\)\cdot \(\overline u^\eps_m \eps\nabla u^\eps_m +
  u^\eps_m \eps\nabla\overline u^\eps_m \)\\
&= \eps\int
  \nabla \(\DIV v\)\cdot\( \overline u^\eps_m (\eps\nabla-iv)u^\eps_m
  +u^\eps_m \overline {(\eps\nabla-iv)u^\eps}_m \)\\
& =\O\( K^\eps_m +\eps^2\), 
\end{align*}
where we have used the conservation of mass and Young's
inequality. From now on, we use the convention that the constant
associated to the notation $\O$ is independent of $m$ 
and $\eps$. We have written the time derivative of $K^\eps_m$ as the
sum of nine terms. The first one corresponds to the conservation of
the energy, and will be canceled by the first term of the time
derivative of $H^\eps_m-K^\eps_m$. We have just bounded the eighth
one. We consider next the sum of four of the seven remaining terms:
the second, third, fifth and seventh,
\begin{align*}
  &\int\big( -\frac{1}{2} |v|^2 \DIV J^\eps_m - \sum_{j,k} \rho^\eps_m v_j v_k
  \d_j v_k +  (v\cdot \nabla v)\cdot J^\eps_m 
-\sum_{j,k}  
  \d_k v_j \RE \( \eps \d_j \overline 
  u^\eps_m \eps \d_k u^\eps_m \) \big)\\
&= \sum_{j,k}\int \Big( v_k \d_j v_k
  \(J^\eps_m\)_j - \lvert u_m^\eps\rvert^2 v_jv_k \d_jv_k + v_j
  \d_jv_k  \(J^\eps_m\)_k- \d_j v_k \RE \( \eps \d_k \overline 
  u^\eps_m \eps \d_j u^\eps_m \).
\end{align*}
Factoring out the term $\d_j v_k$, and recalling that 
\begin{equation*}
 J^\eps_m =\IM\(\overline u_m^\eps \eps \nabla u_m^\eps\), 
\end{equation*}
the above sum can be simplified to:
\begin{equation*}
   -\sum_{j,k}\int \d_j v_k \RE\( \eps\d_j u_m^\eps -iv_j
    u_m^\eps\)\(\overline{\eps\d_k u_m^\eps -iv_k
    u_m^\eps}\)= \O\( K^\eps_m\), 
\end{equation*}
from Cauchy--Schwarz inequality, since $\nabla v\in
L^\infty([0,T]\times\R^n)$. 
We are now left with:
\begin{equation*}
   \frac{d}{dt}K^\eps_m =  \O\( K^\eps_m +\eps^2\) -\frac{d}{dt} \int
  F_m(\rho^\eps_m)+  \int \nabla 
  f_m(\rho^\eps_m) \cdot(\rho^\eps_m v) 
- \int \nabla f(\rho)\cdot \(
  \rho^\eps_m v - 
  J^\eps_m\).  
  \end{equation*}
Since $G'_m(y)=yf'_m(y)$, we infer:
\begin{equation*}
   \frac{d}{dt}K^\eps_m = \O\( K^\eps_m +\eps^2\) -\frac{d}{dt} \int
  F_m(\rho^\eps_m)-  \int 
  G_m(\rho^\eps_m) \DIV v 
- \int \nabla f(\rho)\cdot \(
  \rho^\eps_m v - 
  J^\eps_m\).  
  \end{equation*}
Direct computations yield
\begin{align*}
  \frac{d}{dt} \(H^\eps_m-K^\eps_m\) =&\frac{d}{dt} \int
  F_m\(\rho^\eps_m\) -\int \nabla
  f_m(\rho)\cdot J_m^\eps + \int \( \rho^\eps_m
  - \rho\) v\cdot \nabla f_m\(\rho\)\\
& +\int\( \rho^\eps_m
  - \rho\) G'_m \(\rho\)\DIV v. 
\end{align*}
We therefore come up with:
\begin{equation*}
  \begin{aligned}
    \frac{d}{dt}H^\eps_m=& \O\( K^\eps_m +\eps^2\)  -\int \(
    G_m(\rho^\eps_m)-  G_m(\rho)
    -(\rho^\eps_m -\rho)G_m'(\rho)\)\DIV v\\
&+ \int \nabla\( f(\rho)-f_m(\rho)\)\cdot (J^\eps_m - \rho^\eps_m v).  
  \end{aligned}
\end{equation*}
Note that $f(\rho)-f_m(\rho)\to 0$ in $L^\infty([0,T];W^{1,\infty})$
as $m\to \infty$. We can thus write:
\begin{equation}\label{eq:modul19}
  \begin{aligned}
        \frac{d}{dt}H^\eps_m= &\O\( K^\eps_m +\eps^2\)+o_{m\to
 \infty}(1)\\
 & -\int \(
    G_m(\rho^\eps_m)-  G_m(\rho)
    -(\rho^\eps_m -\rho)G_m'(\rho)\)\DIV v.
  \end{aligned}
\end{equation}
We conclude thanks to the following lemma, whose proof is postponed to the
 end of this section:
\begin{lemma}\label{lem:convex}
There exists $K>0$ independent of $m$ such that  $\forall \rho',\rho \ge 0$,
\begin{equation*}
\left\lvert G_m(\rho')-  G_m(\rho)
    -(\rho' -\rho)G_m'(\rho)\right\rvert \le K \left\lvert F_m(\rho')-
    F_m(\rho) -(\rho' -\rho)F_m'(\rho)\right\rvert .
\end{equation*}
\end {lemma}
Using this lemma and \eqref{eq:modul19}, we infer:
\begin{equation*}
 \frac{d}{dt}H^\eps_m\le  C \(H^\eps_m +\eps^2\)+o_{m\to
 \infty}(1),
\end{equation*}
for some $C$ independent of $m$.
Gronwall lemma yields
\begin{equation*}
  \sup_{t\in [0,T]} H^\eps_m(t) \le C'\eps^2 +o_{m\to \infty }(1),
\end{equation*}
for some $C'$ independent of $m$.
Letting $m\to \infty$, Fatou's lemma yields 
\begin{equation*}
  \sup_{t\in [0,T]} H^\eps(t)\le C' \eps^2.
\end{equation*}
Theorem~\ref{theo:og} then follows from \eqref{eq:convexfinal}.  
\end{proof}

\begin{proof}[Proof of Lemma~\ref{lem:convex}] Taylor's formula yields
\begin{align*}
G_m(\rho')-  G_m(\rho)
    -(\rho' -\rho)G_m'(\rho) &= \( \rho' -\rho\)^2 \int_0^1
    (1-\theta)G_m''\(\rho + \theta(\rho' -\rho)\)d\theta.\\ 
    F_m(\rho')-  F_m(\rho)
    -(\rho' -\rho)F_m'(\rho) &= \( \rho' -\rho\)^2 \int_0^1
    (1-\theta)F_m''\(\rho + \theta(\rho' -\rho)\)d\theta. 
\end{align*}
By definition,
\begin{equation*}
F''_m(y)= f'_m(y)\quad ; \quad G''_m(y) = f'_m (y)+ y f''_m(y).
\end{equation*}
Set, for $y\ge 0$, $h(y)=y^\si /(1+y^\si)$:
\begin{equation*}
 f'_m(y)=\delta_m^{1-\si} h'(\delta_m y)\quad ; \quad
 f''_m(y)=\delta_m^{2-\si} h''(\delta_m y). 
\end{equation*}
Moreover,
\begin{equation*}
h'(y)= \frac{\si y^{\si-1}}{(1+y^\si)^2}\ge 0.
\end{equation*}
Therefore, to prove the lemma, it suffices to show that for all $y\ge 0$,
\begin{equation}\label{eq:23h25}
\left\lvert y h''(y)\right\rvert \le C  h'(y).
\end{equation}
We check the identity
\begin{equation*}
y h''(y) = h'(y )\times\frac{\si-1 -(\si+1)y^\si}{1+y^\si}.
\end{equation*}
The estimate \eqref{eq:23h25} is then straightforward, and
the lemma follows. 
\end{proof}

\section{End of the proof of Theorem~\ref{theo:new}}
\label{sec:concl}

To conclude, the heuristic argument is as follows. From
Theorem~\ref{theo:og}, we expect
\begin{equation*}
  \left\| \eps \nabla u^\eps(t)\right\|_{L^2} \thickapprox \left\|
  v(t) u^\eps(t)\right\|_{L^2} \thickapprox \left\|
  v(t) a(t)\right\|_{L^2}.
\end{equation*}
This follows easily from H\"older's inequality. For the values
$k\in]0,1[$ in Theorem~\ref{theo:og}, we morally use an estimate of
the form
\begin{equation*}
  \left\| |v(t)|^k u^\eps(t)\right\|_{L^2}\lesssim 
\left\| |\eps D_x|^k u^\eps(t)\right\|_{L^2}+ 
\left\| |\eps D_x -v(t)|^k u^\eps(t)\right\|_{L^2},
\end{equation*}
where the first term of the right-hand side goes to zero by
interpolation between  $k=0$ and $k=1$. The aim of the
following lemma is to justify such a  statement. 
\begin{lemma}\label{prop:micro}There exists a constant $K$ such that, for all 
$\eps\in ]0,1]$, for all $s\in [0,1]$, 
for all $u\in H^{1}(\R^n)$ and for all $v\in W^{1,\infty}(\R^n)$,
\begin{equation*}
\| |v|^s u\|_{L^2}\le \| |\eps D_{x}|^s u\|_{L^2}+\| (\eps\nabla -i
v)u\|_{L^2}^{s} 
\| u\|_{L^2}^{1-s}+
\eps^{s/2} K \(1+\| \nabla v\|_{L^\infty}\)\| u\|_{L^2}.
\end{equation*}
\end{lemma}
\begin{proof} 
We begin with the following elementary inequality: 
For all $(x,y)\in \R^{n}\times\R^{n}$ and all $s\in [0,1]$, there holds
\begin{equation}\label{elementary}
| x |^s \le | y|^s + |x-y|^s.
\end{equation}
To see this, note that the result is obvious if $|x|\le |y|$. 
Else, write $|y|=\lambda |x|$ with $\lambda\in [0,1]$ and use the
inequalities $\lambda\le \lambda^s$ and $(1-\lambda)\le (1-\lambda)^s$.  

With this preliminary established, 
introduce the wave-packets operator (see e.g. \cite{CF,Delort,Martinez})
$$
W^\varepsilon v (x,\xi) = c_n \varepsilon^{-3n/4} 
\int _{\R^n} e^{i(x-y)\cdot\xi/\eps  - (x-y)^2/2\eps} v (y) \, dy,
$$
with $c_n=2^{-n/2} \pi^{-3n/4}$. 
The mapping $v\mapsto W^\eps v$ is continuous from the Schwartz class
$\mathcal{S}(\R^n)$  
to $\mathcal{S}(\R^{2n})$, 
and $W^\eps$ extends as an isometry from 
$L^2(\R^n)$ to $L^2(\R^{2n})$:
$$
\| W^\eps v\|_{L^2(\R^{2n})}=\| v\|_{L^2(\R^{n})}.
$$
By applying \eqref{elementary}, we have
$$
\big\| |v(x)|^s W^\eps u\big\|_{L^2(\R^{2n})}\le \big\| |\xi |^s
W^\eps u \big\|_{L^2(\R^{2n})}  
+ \big\| | \xi - v(x)|^s W^\eps u\big\|_{L^2(\R^{2n})}.
$$
Therefore, since
\begin{align*}
 \| | \xi - v(x)|^s W^\eps u\|_{L^2(\R^{2n})} &\le 
 \|  W^\eps u\|_{L^2(\R^{2n})}^{1-s} 
 \big\| | \xi - v(x)| W^\eps u\big\|_{L^2(\R^{2n})}^{s}\\
& \le  \|  u\|_{L^2(\R^{n})}^{1-s} 
 \big\| ( \xi - v(x)) W^\eps u\big\|_{L^2(\R^{2n})}^{s},
\end{align*}
to obtain the desired estimate, we need only prove:
\begin{align}
&\big\| |v(x)|^s W^\eps u- W^\eps(|v|^su) \big\|_{L^2(\R^{2n})} \le
K\eps^{s/2} \| \nabla v\|_{L^\infty}^{s}\| u\|_{L^2},\label{last1}\\ 
&\big\| |\xi|^s W^\eps u - W^\eps(|\eps D_x|^s u)
\big\|_{L^2(\R^{2n})} \le K \eps^{s/2} \| u\|_{L^2},\label{last2}\\ 
& \big\| (i\xi - iv) W^\eps u - W^\eps \bigl(( \eps\nabla -
iv)u\bigr) \big\|_{L^2(\R^{2n})} 
\le K\eps^{1/2} (1+ \| \nabla v\|_{L^\infty})\| u\|_{L^2}.\label{last3}
\end{align}
These properties follows from the fact that the wave packets operator
conjugates the action of  
pseudo-differential operators, approximately, to multiplication by
symbols. For smooth symbols, one has  
sharp results (see \cite{CF,Delort,Martinez}). For the rough symbols
$|v(x)|^s$ and $|\eps\xi|^s$, one can proceed as follows.

To prove~\eqref{last1}, directly from the definition, we compute
\begin{align*}
&\left\| |v|^s W ^\varepsilon u - W ^ \varepsilon (|v|^s u)
\right\|^2_{L^2 (\R^{2n})} \\
&\qquad
=c_{n}^2(2\pi)^n \varepsilon ^{-n/2} \iint  e^{-(x-y)^2
  /\varepsilon } \left\lvert |v(x)|^s  - |v(y)|^s\right\rvert ^2 \left| 
u(y) \right| ^2 \, dy dx.
\end{align*}
Consequently, since $v\in W^{1,\infty}(\R^n)$, the inequality
\eqref{elementary} implies
\begin{align*}
&\left\| |v|^s W ^\varepsilon u - W ^ \varepsilon (|v|^s u)
\right\|^2_{L^2 (\R^{2n})} \\
&\qquad\le K\| \nabla v\|_{L^\infty}^{2s}
\iint \varepsilon ^{-n/2} e^{-(x-y)^2 /\varepsilon } \left|
  x-y\right|^{2s} 
\left| u(y)\right|^2 \, dy dx\\
&\qquad\le K \| \nabla v\|_{L^\infty}^{2s} \iint e^{-z^2} \left|
  \sqrt{\eps}z\right|^{2s} 
\left| u(x-\sqrt{\eps}z) \right|^2 \,dz dx,
\end{align*}
which proves~\eqref{last1}.
We next compute 
$W^\eps(|\eps D_x|^s u) (x,\xi)$: it is given by
\begin{align*}
c_n (2\pi)^{-n/2} \eps^{-7n/4}\iint e^{i(x-y)\cdot
  (\xi-\theta)/\eps-(x-y)^2/2\eps}e^{ix\cdot\theta/\eps} 
|\theta|^s\widehat{u}\(\frac{\theta}{\eps}\)\,d\theta dy,
\end{align*}
where $\widehat{u}$ is the Fourier transform of $u$. Hence, by using 
$$ 
(2\pi)^{-n/2}\int e^{i(x-y)\cdot (\xi-\theta)/\eps-(x-y)^2/2\eps}\, dy =
\eps^{n/2}e^{-(\xi-\theta)^2/2\eps}, 
$$ 
we find
$$
W^\eps\bigl( |\eps D_x|^s u\bigr) (x,\xi)\defn e^{i x\cdot\xi /\eps} 
W^\eps w^\eps (\xi,-x),
$$
with $w^\eps (\tau)\defn | \tau|^s \eps^{-n/2}\widehat{u}(\tau/\eps)$. 
This leads us back to the situation of the previous step (with
$|v(x)|^s$ replaced with $|x|^s$),  and hence \eqref{last2} is proved.  
\smallbreak

Finally, the arguments establishing \eqref{last1} and
\eqref{last2} also yield the usual estimates 
\begin{align*}
&\big\| v W^\eps u - W^\eps(v u ) \big\|_{L^2(\R^{2n})} \le
K\eps^{1/2}\| \nabla v\|_{L^\infty}\| u\|_{L^2},\\ 
&\big\| i\xi W^\eps u - W^\eps(\eps \nabla u)
\big\|_{L^2(\R^{2n})} \le K \eps^{1/2} \| u\|_{L^2}, 
\end{align*}
which proves \eqref{last3}. This completes the proof of the lemma.
\end{proof}
We infer that the heuristic argument of the beginning of this section
is justified:
\begin{corollary}For all $t\in [0,T]$ and all $k\in]0,1]$, we have:
\begin{equation}\label{claimEP}
\liminf_{\eps\to 0} \left\| |\eps D_x|^k u^\eps (t)\right\|_{L^2} \ge 
\left\| |v(t)|^k a(t) \right\|_{L^2}.
\end{equation}
\end{corollary}
\begin{proof} 
Let $t\in [0,T]$. It follows from the previous lemma that
$$
\bigl\| |\eps D_x|^k u^\eps (t)\bigr\|_{L^2} = \| |v(t)|^k u^\eps(t)
\|_{L^2}+o(1). 
$$
Write
$$
\left\| |v(t)|^{k}a(t)\right\|_{L^2} \le \left\| |v(t)|^{k}
  u^\eps(t)\right\|_{L^2}+ 
\left\|  |v(t)|^{2k} \(|u^\eps(t)|^2- |a(t)|^2\)\right\|_{L^1}. 
$$
From H\"older's inequality, the last term is bounded by
\begin{equation}\label{eq:14h27}
\left\| |v(t)|^{2k}\right\|_{L^{1+1/\si}} \left\| |u^\eps(t)|^2-
  |a(t)|^2\right\|_{L^{\si+1}}.   
\end{equation}
When $k\ge \si/(\si+1)$, Lemma~\ref{lem:muk} and Sobolev embedding show that
the first term is bounded on $[0,T]$. 
When $0<k< \si/(\si+1)$, H\"older's inequality yields:
\begin{equation*}
 \left\| |v(t)|^{2k}\right\|_{L^{1+1/\si}}\le  C_N \left\|
 \<x\>^{N}v(t)\right\|_{L^2}^{2k\si/(\si +1)}\quad \text{for
 }N>\frac{n}{2k}\(\frac{\si}{\si+1}-k\). 
\end{equation*}
Lemma~\ref{lem:muk} and Theorem~\ref{theo:og} show
that \eqref{eq:14h27} goes to zero as $\eps$ tends to $0$. 
\end{proof}
To complete the proof of Theorem~\ref{theo:reduc}, 
it remains only to prove that the right-hand side in~\eqref{claimEP}
is non trivial. To see this, we note that, 
from \eqref{eq:limitev},
\begin{equation*}
 a_{\mid t=0}=a_0\quad ;\quad v_{\mid t=0} =0\quad ;\quad \d_t v_{\mid
 t=0} = -\nabla \( |a_0|^{2\si}\). 
\end{equation*}
Therefore, by continuity (see Lemma~\ref{lem:muk}), we obtain the
following result. 
\begin{lemma}
There exists $\tau>0$ such that
\begin{equation}\label{eq:decolle}
  \int |v(\tau,x)|^{2k} |a(\tau,x)|^2dx>0, \quad \forall k\in [0,1]. 
\end{equation}
\end{lemma}

This implies Theorem~\ref{theo:reduc}, hence
Theorem~\ref{theo:new}. 

\begin{remark}\label{rem:semiquasi}
  We can compare the results of this paper with the analysis in
  \cite{CCT2}. The approximate solution used in \cite{CCT2} consists
  in neglecting the Laplacian in \eqref{eq:nlssemi}:
  \begin{equation*}
    i\eps\d_t w^\eps = |w^\eps|^{2\si}w^\eps \quad ;\quad w^\eps_{\mid
    t=0}=a_0,\quad \text{hence}\quad w^\eps(t,x)=a_0(x)e^{-it
    |a_0(x)|^{2\si}/\eps}.   
  \end{equation*}
A direct application of Gronwall lemma shows that $w^\eps$
  is  a suitable 
  approximation of $u^\eps$ up to time of order $c\eps |\log
  \eps|^\theta$,  for some $c,\theta>0$. The Taylor expansion in time for $v$
  shows that
  \begin{equation*}
    v(t,x) = -t \nabla \(|a_0(x)|^{2\si}\) +\O\(t^3\).
  \end{equation*}
The formal analysis of \cite[\S 3.1]{CaARMA} is thus
justified also in this case: $w^\eps(t)$ is a good approximation of
$u^\eps(t)$ for $t\ll \eps^{1/3}$:
  \begin{equation*}
    \||\eps D_x|^s u^\eps(t)\|_{L^2}\thickapprox \||v(t)|^s
    a(t)\|_{L^2}\thickapprox  \||\eps D_x|^s w^\eps(t)\|_{L^2}\quad
    \text{for }t\ll \eps^{1/3}. 
  \end{equation*}
 To prove this point, it seems
  necessary to perform a quasilinear analysis (see~\S\ref{sec:muk}),
  and the semilinear approach based on Gronwall lemma is not enough.  
\end{remark}
\section{Final remarks}\label{sec:geometrie}

To conclude this paper, we note that the approach presented here
remains efficient in  the case of non-trivial geometries. Indeed, the
scaling argument that we have used in \S\ref{sec:reduc} is merely
helpful for the intuition, to guess a suitable approximate
solution. This argument meets the strategy adopted in the appendix of
\cite{BGTENS}. Introduce $\om^h$ and $A^h$ given by
\begin{equation*}
  v(t,x) = \om^h\(h^2\eps t,hx\)\quad ;\quad a(t,x) =
  h^{\frac{n}{2}-s}A^h\(h^2\eps t,hx\). 
\end{equation*}
The key approximation that we have used,
\begin{equation*}
  \left\lvert \eps\nabla u^\eps(t,x)\right\rvert^2 \approx \left\lvert
  v(t,x)a(t,x) \right\rvert^2,\quad 0\le t\le T,
\end{equation*}
then reads
\begin{equation*}
  \left\lvert \eps h \nabla \psi^h (t,x)\right\rvert^2 \approx \left\lvert
  \om^h(t,x)A^h(t,x) \right\rvert^2, \quad 0\le t\le h^2 \eps T.
\end{equation*}
More precisely, in terms of the initial problem \eqref{eq:nls},
Theorem~\ref{theo:og} reads:
\begin{equation*}
  \begin{aligned}
   & h^{-s} \left\lVert \(\eps h \nabla
    -i\om^h\)\psi^h\right\rVert_{L^\infty([0,h^2\eps T];L^2)}^2 \\
+&h^{\(\frac{n}{2}-s\)2(\si+1) -\frac{n}{2}}\left\lVert \(\lvert
    \psi^h\rvert^2 - \lvert A^h\rvert^2\)^2\(\lvert
    \psi^h\rvert^{2\si-2} +\lvert
    A^h\rvert^{2\si-2}\)\right\rVert_{L^\infty([0,h^2\eps
    T];L^1)}\\
&=\O\(\eps^2\). 
  \end{aligned}
\end{equation*}
It is essentially this estimate that we have used to prove
Theorem~\ref{theo:new} (and Corollary~\ref{cor:energy} stems exactly
from this estimate). 
\smallbreak

Suppose for instance that $x$ belongs to a bounded
domain $M$, and not to all of $\R^n$, and that
we consider 
\eqref{eq:nls}  on $M$, with Dirichlet or Neumann boundary
condition. If $a_0$ is compactly supported in a ball, 
contained in the interior
of $M$, 
\begin{equation*}
  \operatorname{supp} a_0\subset B \Subset M,
\end{equation*}
then we can still
consider \eqref{eq:limite}, viewed as a 
system on $\R^n$. The key remark is that for $a_0$ compactly
supported, the smooth solutions to \eqref{eq:limite} have a finite
speed of propagation, which is \emph{zero}. This is an important step
in proving that smooth solutions develop singularities in finite time;
see \cite{MUK86,Xin98}. Therefore, $(\phi,a)$  remains supported in
$B$ for $t\in [0,T]$; up to changing the origin to the center of $B$,
so does $(\om^h,A^h)$ on the time interval $[0,h^2\eps T]$. In
particular, $(\om^h,A^h)$ is supported away from the boundary of $M$
for $t\in [0,h^2\eps T]$.
\smallbreak

In the proof of Theorem~\ref{theo:og}, 
the integrations by parts affect $v$ or $a$, that is, $\om^h$ or
$A^h$. Only two terms 
in the differentiation of the modulated energy functional do not
contain $\om^h$ or $A^h$: these two terms correspond to the global
energy of $\psi^h$, which is a non-increasing function of time for
strong solutions in the case of $\R^n$ (since it is in fact
constant). As a matter of fact, this property is needed for 
strong solutions only, since it remains for weak solutions, by Fatou's
lemma. Therefore, Theorem~\ref{theo:new} remains valid
on $M$, provided that we can construct strong solutions with a
non-increasing energy. This is the case of compact surfaces  when
$\si\ge 1$, and of compact three dimensional manifolds 
when $\si=1$, see \cite{BGT}. This is also the case of bounded domains
in $\R^2$ for $\si\ge 1$, see \cite{RamonaBSMF}, of the ball in $\R^3$
for $\si= 1$ and radial data \cite{RamonaCPDE}, and of exterior
domains in $\R^3$, 
for $\si=1$ \cite{RamonaJMPA}. 
\smallbreak

Note however that the notion of criticality
may differ on a curved space (see e.g. \cite{BGTMRL,ThomannSurf}): the
curvature of a manifold may create more ill-posedness phenomena, but
since in the proof of ill-posedness in \cite{CCT2} (see also
\cite{CaARMA}), the Laplacian is 
neglected, the critical Sobolev exponent for local well-posedness cannot
be less than in the case of $\R^n$. 
\smallbreak

This approach suggests that we can consider a
more general manifold, up to working on a local chart, and provided
that the energy associated to strong solutions of \eqref{eq:nls} is
a non-increasing function of time, a question which we leave out at
this stage. 
\smallbreak

Finally, we go back to the whole space case, $x\in \R^n$. As recalled
above, if $a_0$ is compactly  supported, then the ansatz that we
consider remains supported in the same compact so long as the solution
to \eqref{eq:limiteS} remains smooth; see \cite{MUK86}, and also
\cite{Xin98}. The justification of WKB analysis for short time shows
that at least when $\si=1$ (\cite{Grenier98}), or $\si \in \N$ and $n\le 3$
(\cite{AC-BKW}), we have, thanks to Borel lemma,
\begin{equation*}
  u^\eps(t,x) = u^\eps_{\rm app}(t,x) + \O\(\eps^\infty\), \quad
  \text{in }C\([0,T];L^2\cap L^\infty\), 
\end{equation*}
where $u^\eps_{\rm app}$ is supported in the same compact as
$a_0$. This seems to be an encouraging remark, in view of considering
\emph{fixed} initial data as in \cite{Lebeau05}, instead of a sequence
of initial data like here. However, the information that we do not
have for the Schr\"odinger equation, and which is available for the
wave equation, is a notion of finite speed of propagation for
\emph{weak} solutions to the nonlinear equation. This seems to be the
only obstacle to consider fixed initial data in the Schr\"odinger
case.

\subsection*{Acknowledgments} The authors are grateful to
Patrick G\'erard for stimulating comments on this work. 

\bibliographystyle{amsplain}
\bibliography{../MathAnnalen/loss}

\providecommand{\bysame}{\leavevmode\hbox to3em{\hrulefill}\thinspace}
\providecommand{\MR}{\relax\ifhmode\unskip\space\fi MR }
\providecommand{\MRhref}[2]{%
  \href{http://www.ams.org/mathscinet-getitem?mr=#1}{#2}
}
\providecommand{\href}[2]{#2}
\begin{thebibliography}{10}

\bibitem{AC-BKW}
T.~Alazard and R.~Carles, \emph{Supercritical geometric optics for nonlinear
  {S}chr\"odinger equations}, archived as {\tt arXiv:0704.2488}, 2007.

\bibitem{AlinhacGerardUS}
S.~Alinhac and P.~G\'erard, \emph{Pseudo-differential operators and the
  {N}ash-{M}oser theorem}, Graduate Studies in Mathematics, vol.~82, American
  Mathematical Society, Providence, RI, 2007, Translated from the 1991 French
  original by Stephen S. Wilson.

\bibitem{RamonaCPDE}
R.~Anton, \emph{Cubic nonlinear {S}chr\"odinger equation on three dimensional
  balls with radial data}, Comm. Partial Differential Equations, to appear.
  Archived as {\tt arXiv:math/0608689}.

\bibitem{RamonaJMPA}
\bysame, \emph{Global existence for defocusing cubic {NLS} and
  {G}ross-{P}itaevskii equations in three dimensional exterior domains}, J.
  Math. Pures Appl. (9), to appear. Archived as {\tt arXiv:math.AP/0701304}.

\bibitem{RamonaBSMF}
\bysame, \emph{Strichartz inequalities for {L}ipschitz metrics on manifolds and
  nonlinear {S}chr\"odinger equation on domains}, Bull. Soc. Math. France, to
  appear. Archived as {\tt arXiv:math.AP/0512639}.

\bibitem{BrenierCPDE}
Y.~Brenier, \emph{Convergence of the {V}lasov-{P}oisson system to the
  incompressible {E}uler equations}, Comm. Partial Differential Equations
  \textbf{25} (2000), no.~3-4, 737--754.

\bibitem{BGTMRL}
N.~Burq, P.~G\'erard, and N.~Tzvetkov, \emph{An instability property of the
  nonlinear {S}chr\"odinger equation on {$S\sp d$}}, Math. Res. Lett.
  \textbf{9} (2002), no.~2-3, 323--335.

\bibitem{BGT}
\bysame, \emph{Strichartz inequalities and the nonlinear {S}chr\"odinger
  equation on compact manifolds}, Amer. J. Math. \textbf{126} (2004), no.~3,
  569--605.

\bibitem{BGTENS}
\bysame, \emph{Multilinear eigenfunction estimates and global existence for the
  three dimensional nonlinear {S}chr\"odinger equations}, Ann. Sci. \'Ecole
  Norm. Sup. (4) \textbf{38} (2005), no.~2, 255--301.

\bibitem{CaARMA}
R.~Carles, \emph{Geometric optics and instability for semi-classical
  {S}chr\"odinger equations}, Arch. Ration. Mech. Anal. \textbf{183} (2007),
  no.~3, 525--553.

\bibitem{CaBKW}
\bysame, \emph{{WKB} analysis for nonlinear {S}chr\"odinger equations with
  potential}, Comm. Math. Phys. \textbf{269} (2007), no.~1, 195--221.

\bibitem{CW90}
T.~Cazenave and F.~Weissler, \emph{The {C}auchy problem for the critical
  nonlinear {S}chr\"odinger equation in ${H}^s$}, Nonlinear Anal. TMA
  \textbf{14} (1990), 807--836.

\bibitem{CG05}
C.~Cheverry and O.~Gu\`es, \emph{Counter-examples to concentration-cancellation
  and supercritical nonlinear geometric optics for the incompressible {E}uler
  equations}, Arch. Ration. Mech. Anal., to appear.

\bibitem{CCT2}
M.~Christ, J.~Colliander, and T.~Tao, \emph{Ill-posedness for nonlinear
  {S}chr\"odinger and wave equations}, archived as {\tt arXiv:math.AP/0311048}.

\bibitem{CKSTTAnnals}
J.~Colliander, M.~Keel, G.~Staffilani, H.~Takaoka, and T.~Tao, \emph{Global
  well-posedness and scattering for the energy--critical nonlinear
  {S}chr\"odinger equation in {$\mathbb R\sp 3$}}, Ann. of Math. (2), to
  appear.

\bibitem{CF}
A.~C{\'o}rdoba and C.~Fefferman, \emph{Wave packets and {F}ourier integral
  operators}, Comm. Partial Differential Equations \textbf{3} (1978), no.~11,
  979--1005.

\bibitem{Delort}
J.-M. Delort, \emph{F.{B}.{I}. transformation. second microlocalization and
  semilinear caustics}, Lecture Notes in Mathematics, vol. 1522,
  Springer-Verlag, Berlin, 1992.

\bibitem{PGX93}
P.~G{\'e}rard, \emph{Remarques sur l'analyse semi-classique de l'\'equation de
  {S}chr\"odinger non lin\'eaire}, S\'eminaire sur les \'Equations aux
  D\'eriv\'ees Partielles, 1992--1993, \'Ecole Polytech., Palaiseau, 1993,
  pp.~Exp.\ No.\ XIII, 13.

\bibitem{GV84}
J.~Ginibre and G.~Velo, \emph{On the global {C}auchy problem for some nonlinear
  {S}chr\"odinger equations}, Ann. Inst. H. Poincar\'e Anal. Non Lin\'eaire
  \textbf{1} (1984), no.~4, 309--323.

\bibitem{GV85c}
\bysame, \emph{The global {C}auchy problem for the nonlinear {S}chr\"odinger
  equation revisited}, Ann. Inst. H. Poincar\'e Anal. Non Lin\'eaire \textbf{2}
  (1985), 309--327.

\bibitem{Grenier98}
E.~Grenier, \emph{Semiclassical limit of the nonlinear {S}chr\"odinger equation
  in small time}, Proc. Amer. Math. Soc. \textbf{126} (1998), no.~2, 523--530.

\bibitem{Kato87}
T.~Kato, \emph{On nonlinear {S}chr\"odinger equations}, Ann. IHP (Phys.
  Th\'eor.) \textbf{46} (1987), no.~1, 113--129.

\bibitem{KatoNLS}
\bysame, \emph{Nonlinear {S}chr\"odinger equations}, Schr\"odinger operators
  (S\o nderborg, 1988), Lecture Notes in Phys., vol. 345, Springer, Berlin,
  1989, pp.~218--263.

\bibitem{Lebeau01}
G.~Lebeau, \emph{Non linear optic and supercritical wave equation}, Bull. Soc.
  Roy. Sci. Li\`ege \textbf{70} (2001), no.~4-6, 267--306 (2002), Hommage \`a
  Pascal Laubin.

\bibitem{Lebeau05}
\bysame, \emph{Perte de r\'egularit\'e pour les \'equations d'ondes
  sur-critiques}, Bull. Soc. Math. France \textbf{133} (2005), 145--157.

\bibitem{LinZhang}
F.~Lin and P.~Zhang, \emph{{S}emiclassical limit of the {G}ross-{P}itaevskii
  equation in an exterior domain}, Arch. Rational Mech. Anal. \textbf{179}
  (2006), no.~1, 79--107.

\bibitem{MUK86}
T.~Makino, S.~Ukai, and S.~Kawashima, \emph{Sur la solution \`a support compact
  de l'\'equation d'{E}uler compressible}, Japan J. Appl. Math. \textbf{3}
  (1986), no.~2, 249--257.

\bibitem{Martinez}
A.~Martinez, \emph{An introduction to semiclassical and microlocal analysis},
  Universitext, Springer-Verlag, New York, 2002.

\bibitem{RV}
E.~Ryckman and M.~Visan, \emph{Global well-posedness and scattering for the
  defocusing energy--critical nonlinear {S}chr\"odinger equation in {$\mathbb
  R\sp{1+4}$}}, Amer. J. Math. \textbf{129} (2007), no.~1, 1--60.

\bibitem{TaoDisp}
T.~Tao, \emph{Nonlinear dispersive equations}, CBMS Regional Conference Series
  in Mathematics, vol. 106, Published for the Conference Board of the
  Mathematical Sciences, Washington, DC, 2006, Local and global analysis.

\bibitem{Taylor3}
M.~Taylor, \emph{Partial differential equations. {III}}, Applied Mathematical
  Sciences, vol. 117, Springer-Verlag, New York, 1997, Nonlinear equations.

\bibitem{ThomannAnalytic}
L.~Thomann, \emph{Instabilities for supercritical {S}chr\"odinger equations in
  analytic manifolds}, J. Differential Equations, to appear. Archived as {\tt
  arXiv:0707.1785}.

\bibitem{ThomannSurf}
\bysame, \emph{The {WKB} method and geometric instability for non linear
  {S}chr\"odinger equations on surfaces}, Bull. Soc. Math. France, to appear.
  Archived as {\tt arXiv:math/0609805}.

\bibitem{Xin98}
Z.~Xin, \emph{Blowup of smooth solutions of the compressible {N}avier-{S}tokes
  equation with compact density}, Comm. Pure Appl. Math. \textbf{51} (1998),
  229--240.

\end{thebibliography}

\end{document}